\documentclass[a4paper, 12pt]{article}
\usepackage[utf8]{inputenc}%keine Ahnung
\usepackage{a4}%kleinere Ränder
\usepackage[T1]{fontenc}%Silbentrennung
\usepackage[english]{babel}

\usepackage{amssymb,amsmath,amsthm}
\usepackage{enumerate}
\usepackage{graphicx}
\usepackage{epsfig}
\usepackage{comment}
\usepackage{xspace}
\usepackage{float}
\usepackage{multirow}
\usepackage{extarrows}
\usepackage{mathtools}
\usepackage{color}
\usepackage{nomencl}%Notationsverzeichnis
\usepackage{nicefrac}
\usepackage{cite}
\makenomenclature
\newcommand{\N}{\ensuremath{\mathbb{N}}\xspace}

\newcommand{\R}{\ensuremath{\mathbb{R}}\xspace}

\newcommand{\eps}{\epsilon}
\newcommand{\loc}{\mathrm{loc}}

\renewcommand{\epsilon}{\varepsilon}

\newcommand{\palpha}{\partial_i}
\newcommand{\udel}{u_\delta}
\newcommand{\fudel}{\mathbf{u}_\delta}
\newcommand{\fu}{\overline{\mathbf{u}}}
\newcommand{\vdel}{v_\delta}
\newcommand{\Fdel}{F_\delta}
\newcommand{\thdel}{\Theta_\delta}
\newcommand{\psidel}{\omega_\delta}
\newcommand{\phidel}{\varphi_\delta}
\newcommand{\pvd}{\overline{\palpha\vdel}}

\newcommand{\dx}{\, \mathrm{d}x}

\newcommand{\intom}{\intop_{\Omega}}
\newcommand{\intomd}{\intop_{\Omega-D}}
\newcommand{\intb}{\intop_{B_{R}(x_0)}}
\newcommand{\intdb}{\intop_{B_{2R}(x_0)}}
\newcommand{\nabu}{\nabla u}
\newcommand{\navdel}{\nabla v_\delta}
\newcommand{\nabdel}{\nabla u_\delta}

\newcommand{\gr}[1]{(\ref{#1})}

\newcommand{\leb}{\mathcal{L}^2}

\newcommand{\zkz}{\ensuremath{\mathbb{R}^{2\times 2}}}
\newcommand{\ft}[1]{\mathbf{#1}}
\newcommand{\phitil}{\widetilde{\varphi}_\delta}
\newcommand{\DG}{\intom D^2G(\nabla\udel-\vdel)}
\newcommand{\pnu}{\palpha\nabla\udel}

\newcommand{\wthdel}{\widetilde{\Theta}_\delta}
\newcommand{\wudel}{\widetilde{u}_\delta}
\newcommand{\intomo}{\intop_{\Omega_0}}
\newcommand{\intomod}{\intop_{\Omega_0}}
\newcommand{\hthdel}{\widehat{\Theta}_\delta}
\def\Yint#1{\mathchoice
    {\YYint\displaystyle\textstyle{#1}}%
    {\YYint\textstyle\scriptstyle{#1}}%
    {\YYint\scriptstyle\scriptscriptstyle{#1}}%
    {\YYint\scriptscriptstyle\scriptscriptstyle{#1}}%
      \!\int}
\def\YYint#1#2#3{{\setbox0=\hbox{$#1{#2#3}{\int}$}
    \vcenter{\hbox{$#2#3$}}\kern-.5\wd0}}

\def\mint{\Yint -}

\begin{document}

\numberwithin{equation}{section}
%Eigene Theoremumgebungen
\newtheoremstyle{break}{15pt}{15pt}{\itshape}{}{\bfseries}{}{\newline}{}
\theoremstyle{break}
\newtheorem*{Satz*}{Theorem}
\newtheorem*{Rem*}{Remark}
\newtheorem*{Lem*}{Lemma}
\newtheorem{Satz}{Theorem}[section]
\newtheorem{Rem}{Remark}[section]
\newtheorem{Lem}{Lemma}[section]
\newtheorem{Prop}{Proposition}[section]
\newtheorem{Cor}{Corollary}[section]
\theoremstyle{definition}
\newtheorem{Def}[Satz]{Definition}
\newtheorem*{Ass}{General Assumptions}
\parindent2ex

\newenvironment{rightcases}
  {\left.\begin{aligned}}
  {\end{aligned}\right\rbrace}

\begin{center}{\Large \bf Coupled variational problems of linear growth related to the denoising and inpainting of images}
\end{center}

\begin{center}
J. M\"uller
\end{center}

\begin{abstract}
In this note we present some results that were already conjectured in the work \cite{BFW} by Bildhauer, Fuchs and Weickert, where they have investigated analytical aspects of coupled variational models with applications to mathematical imaging. Here we focus on variants of linear growth, which require a treatment in the framework of relaxation theory and convex analysis. Following basic ideas from \cite{BF2} and \cite{BF3}, we establish existence and regularity  of (dual-)solutions.
\end{abstract}


\begin{section}{Introduction}
The paper at hand takes up  ideas from the work \cite{BFW}, where certain variational problems  from the field of image analysis have been studied. The general application background is given by the task to retrieve  a digital grey-scale image, i.e. a real valued function which maps every point (which can in this context be thought of as an infinitely small ``pixel``) of  a plane domain $\Omega\subset\R^2$ to a grey-value between 0 (indicating a black point) and 1 (indicating a white point), from a flawed observation $f:\Omega-D\rightarrow[0,1]$. Here, in our understanding the term ``flawed`` includes the phenomena of a statistical distortion (called ``noise``) as well as the missing of some  parts of the data, which means that $f$ is only defined outside of a subset $D\subset\Omega$ (the ``deficiency set``). A well established approach to the solution of this problem (called ``pure denoising``, if only the first type of data corruption is considered and ``inpainting`` for the second type) consists in minimizing a suitable energy that penalizes fluctuations of the data. In their fundamental work \cite{ROF}, Rudin, Osher and Fatemi proposed to consider the variational problem
\begin{align}\label{ROF}
\intom |\nabla u|+\frac{\lambda}{2}\intop_{\Omega-D} (u-f)^2\dx\rightarrow\min\text{ in }BV(\Omega)
\end{align}
in the space of functions of bounded variation (see \cite{AFP} or \cite{Giu} for details on this function space). $BV$-functions are very well-suited for modeling objects like images, since they are allowed to have jump-type discontinuities which can reflect edges and sharp contours. On the other hand, they are still in some sense regular enough to allow a treatment with analytical methods. However, the model \gr{ROF} has some drawbacks. First of all, the total variation $\intom|\nabu|$ is neither strictly convex nor differentiable in its argument, hence unsuited to a treatment with PDE-methods. Secondly, numerical simulations show that the solutions of \gr{ROF} are frequently afflicted with the so called ''staircaising effect'' (see e.g. \cite{CMM}), which becomes manifest in piecewise constant regions of the minimizing function (resembling a staircase). The first problem has been circumvented in \cite{BF1} by using the concept of convex functions of a measure (see \cite{DT}). More precisely, the quantity $\intom |\nabu|$ is being replaced with $\intom F(\nabu)$ for a smooth strictly convex function $F:\R^2\rightarrow [0,\infty)$ of linear growth which approximates $|\,.\,|$, e.g. $F(\xi)=\sqrt{\eps^2+|\xi|^2}-\eps$. 

For avoiding the staircasing effect, one could raise the order of differentiability, i.e. consider the problem
\begin{align}\label{ROF2}
\intom |\nabla^2 u|+\frac{\lambda}{2}\intop_{\Omega-D} (u-f)^2\dx\rightarrow\min\text{ in }BV^2(\Omega)
\end{align}
where $\nabla^2u$ is the Hessian matrix and $BV^2(\Omega)$ denotes the set of all functions  $u\in L^1(\Omega)$ such that the weak gradient $\nabla u$ is a $BV$-function. The undesirable effects are then shifted to the first derivative of the solution and therefore become less evident. Analytical properties of these models have been studied in \cite{FM} by Fuchs and the author. However, higher order models are computationally more difficult to handle as they lead to partial  differential equations of at least fourth order. That is why in \cite{BFW}, a different approach has been pursued. There, in place of \gr{ROF2}, a coupled problem has been considered which is obtained by introducing a vector-valued variable $v$, serving as a substitute for the gradient of $u$. To be more precise, the idea is to study the problem
\begin{align}\label{ROFc}
\begin{split}
(u,v)\mapsto\alpha\intom |\nabla v|+\beta\intom |\nabla u-v\cdot\mathcal{L}^n|+\frac{\lambda}{2}\intop_{\Omega-D} (u-f)^2\dx\\
\rightarrow\min\text{ in }BV(\Omega)\times BV(\Omega,\R^n).
\end{split}
\end{align}
Minimizing the middle term, the so called ``coupling term``, entails $v\approx \nabla u$ and hence $\nabla v\approx \nabla^2u$. Solutions of \gr{ROFc} serve as an approximation to those of \gr{ROF2} with the advantage, that the associated system of differential equations is merely of second order. Of course there is plenty of different choices of densities other than $|\cdot|$ in the leading as well as in the coupling term of \gr{ROFc} and actually, in \cite{BFW}, various constellations of power and superlinear growth have been considered. Here, however, we focus on the case where both the leading and the coupling term are of linear growth. Regularity properties of minimizers of functionals of this type have been conjectured in Remark 6.4 of \cite{BFW}. After these few introductory words, we continue with an overview of our assumptions and results.

Let $F:\R^{2\times 2}\rightarrow \R$ and $G:\R^2\rightarrow\R$ be strictly convex, satisfying the following
\begin{Ass} Let $\Omega$ be a Lipschitz Domain in $\R^2$, $D$ a measurable subset such that $\Omega-\overline{D}\neq\emptyset$ and $f:\Omega-D\rightarrow\R$  bounded and measurable. We demand at least ($c$ denotes a generic positive constant):
\begin{enumerate}[\hspace{1cm}(F1)]
 \item \label{F1} $F\in C^2(\zkz)$, $F(-p)=F(p)$, $F(0)=0$, $DF(0)=0$, $|DF|\leq c$,
 \item \label{F2}$0<D^2F(p)(q,q)\leq c \dfrac{1}{1+|p|}|q|^2$ for all $p,q\in\zkz$, 
 \item \label{F3}$c_1|p|-c_2\leq F(p)$ for some $c_1>0$, $c_2\in\R$ and all $p\in\zkz$.
\end{enumerate}
\begin{enumerate}[\hspace{1cm}(G1)]
 \item \label{G1}\label{a} $G\in C^2(\R^2)$, $G(-y)=G(y)$, $G(0)=0$, $DG(0)=0$, $|DG|\leq c$,
 \item \label{G2}$0<D^2G(y)(x,x)\leq c \dfrac{1}{1+|y|}|x|^2$ for all $x,y\in\R^2$, 
 \item \label{G3}$c_1|y|-c_2\leq G(y)$ for some $c_1>0$, $c_2\in\R$ and all $y\in\R^2$.
\end{enumerate}
\end{Ass}
We further define
\begin{align}
 V:=W^{1,1}(\Omega)\times W^{1,1}(\Omega,\R^2).
\end{align}
The underlying problem then reads
\begin{align}\label{E}\tag{$\mathcal{P}$}
\begin{split}
 E(u,v):=\alpha\intom F(\nabla v)\dx+\beta\intom G(\nabla u-v)\dx+\intop_{\Omega-D}(u-f)^2\dx\\
\rightarrow\min\text{ in } V,
\end{split}
\end{align}
where $\alpha$ and $\beta$ are  positive parameters which control the balance between the leading term $\intom F(\nabla v)\dx$ and the coupling term  $\intom G(\nabla u-v)\dx$.
\begin{Rem}
It is needless to mention, that by an iteration procedure the coupling method can be used to reduce functionals of any order higher than two to a problem which involves first derivatives only.
\end{Rem}
Of course, we cannot expect solvability of \gr{E} in the non reflexive space $V$ in general and we therefore have to pursue the approach from \cite{BF2} and \cite{BF3}; which means to consider suitably relaxed variants of the above problem. The first method is the relaxation of \gr{E} in the space $BV(\Omega)\times BV(\Omega,\R^2)$, using the concept of convex functions of a measure (see \cite{DT}).  We therefore replace $E$ with the functional
\begin{align*}
\widetilde{E}(u,v)=\alpha\intom F(\nabla v)+\beta\intom G(\nabla u-v\cdot\mathcal{L}^2)+\intop_{\Omega-D}(u-f)^2\dx
\end{align*}
and  look for solutions of
\begin{align}\label{Etilde}\tag{$\widetilde{\mathcal{P}}$}
\widetilde{E}(u,v)\rightarrow\min\text{ in }BV(\Omega)\times BV(\Omega,\R^2),
\end{align}
where for  finite Radon measures $\mu\in\mathcal{M}(\Omega,\R^2)$, $\nu\in\mathcal{M}(\Omega,\R^{2\times 2})$ we declare (cf. formula (1.4) in \cite{DT})
\begin{align*}
 \intom F(\nu):=\intom F(\nu^a)\dx+\intom F^\infty\left(\frac{d\nu^s}{d|\nu^s|}\right)\,\mathrm{d}|\nu^s|,
\end{align*}
and 
\begin{align*}
 \intom G(\mu):=\intom G(\mu^a)\dx+\intom G^\infty\left(\frac{d\mu^s}{d|\mu^s|}\right)\,\mathrm{d}|\mu^s|.
\end{align*}
With $\mu=\mu^a\mathcal{L}^2+\mu^s$ we denote the Lebesgue decomposition of $\mu$  and $F^\infty(p):=\lim\limits_{t\rightarrow\infty}\frac{F(tp)}{t}$ for all $p\in\R^{2\times 2}$ as well as $G^\infty(y):=\lim\limits_{t\rightarrow\infty}\frac{G(ty)}{t}$ for all $y\in\R^2$. Note, that the Lebesgue decomposition of the measure $\nabla u-v\cdot\mathcal{L}^2$ is given by 
\[
\nabla u-v\cdot\mathcal{L}^2=(\nabla^au-v)\cdot\mathcal{L}^2+\nabla^s u,
\]
where for $u\in BV(\Omega)$, we denote by $\nabla^a u$, $\nabla^s u$ the absolutely continuous and the singular part, respectively, of the gradient measure with respect to $\mathcal{L}^2$.
Then as in \cite{BFW}, Theorem 5.1, we can show
\begin{Satz}\label{Thm1.1}
  Under our general assumptions regarding $\Omega$, $D$, $f$, $F$ and $G$ it holds: 
\begin{enumerate}[a)]
\item Problem \gr{Etilde} has at least one solution $(u,v)\in BV(\Omega)\times BV(\Omega,\R^2)$.
\item If $(u,v)$ and $(\widetilde{u},\widetilde{v})$ are two distinct solutions of \gr{Etilde}, then
\begin{align*}
 &u=\widetilde{u}\text{ a.e. on }\Omega-D,\;\nabla^au-v=\nabla^a\widetilde{u}-\widetilde{v}\text{ a.e. on }\Omega\\
&\hspace{2cm}\text{ and }\nabla^av=\nabla^a\widetilde{v}\text{ a.e. on }\Omega.
\end{align*} 
In particular, if $D=\emptyset$, i.e. in the case of pure denoising, the solution of \gr{Etilde} is unique.
\item The set $\mathcal{M}$ of all solutions of \gr{Etilde} coincides with the set of all $L^1(\Omega)\times L^1(\Omega,\R^2)$-cluster points of $E$-minimizing sequences in $W^{1,1}(\Omega)\times W^{1,1}(\Omega,\R^2)$. If $E$ admits a minimizer $(u,v)$ in the Sobolev class $V$, then $\mathcal{M}=\{(u,v)\}$.
\end{enumerate}
\end{Satz}

Another well established approach towards the relaxation of \gr{E} is via convex duality (cf. \cite{BF3},  \cite{FT} or \cite{MT}). Here, as in \cite{BF3}, we pass to the dual formulation via Lagrangians. In order to simplify our notation, we define the linear operator
\begin{align}
 \Lambda:V\rightarrow Y:=L^1(\Omega,\R^2)\times L^1(\Omega,\R^{2\times 2}),\,\mathbf{u}=(u,v)\mapsto (\nabla u-v,\nabla v)
\end{align}
as well as the function
\[
 \mathcal{F}:\R^2\times\R^{2\times 2}\rightarrow\R,\,(y,p)\mapsto \alpha F(p)+\beta G(y).
\]
We observe, that  problem \gr{E} can now be written for short as
\begin{align}
 E(\mathbf{u})=\intom \mathcal{F}(\Lambda \mathbf{u})\dx+\intop_{\Omega-D}(u-f)^2\dx\rightarrow\min
\end{align}
where $\mathbf{u}=(u,v)\in W^{1,1}(\Omega)\times W^{1,1}(\Omega,\R^2)$.
By means of this representation, it is easy to see how to apply the results from \cite{ET}, Remark 4.1 and 4.2 on pp. 60-61 in order to obtain the problem in duality to \gr{E}. First, for $\ft{w}=(u,v)\in V$ and $\ft{y}=(\kappa,\lambda)\in Y^*=L^\infty(\Omega,\R^2)\times L^\infty(\Omega,\R^{2\times 2})$ we define the associated Lagrangian through
\begin{align}\label{lag}
 \ell(\ft{w},\ft{y}):=\intom \Lambda(\ft{w})\odot \ft{y}\dx-\intom \mathcal{F}^*(\ft{y})\dx+\intomd (u-f)^2\dx,
\end{align}
where for $(x,p),(y,q)\in \R^2\times \R^{2\times 2}$ we set
\begin{align*}
 (x,p)\odot(y,q):=x\cdot y+p:q,
\end{align*}
with ''$\cdot$'' and ''$:$'' denoting the standard scalar products on $\R^2$ and $\R^{2\times 2}$, respectively. Furthermore, $\mathcal{F}^*$ is the convex dual to the function $\mathcal{F}$ which, by Remark 4.3 on p. 61 in \cite{ET} can be split into
\begin{align*}
 \mathcal{F}^*(\kappa,\lambda)=\alpha F^*(\lambda)+\beta G^*(\kappa).
\end{align*}
Hence, we may write \gr{lag} as
\begin{align*}
 \ell(\ft{w},\ft{y})=\intom \nabla v:\kappa+(\nabla u-v)\cdot\lambda\dx-\intom \alpha F^*(\lambda)+\beta G^*(\kappa)\dx+\intomd (u-f)^2\dx
\end{align*}
and it holds (see \cite{ET}, p. 56)
\begin{align*}
 E(\ft{w})=\sup_{\ft{y}\in V^*}\ell(\ft{w},\ft{y}).
\end{align*}
The dual functional $R:Y^*\rightarrow [0,\infty]$ is now defined through
\begin{align}\label{R1}
 R(\ft{y}):=\inf_{\ft{w}\in V}\ell(\ft{w},\ft{y})
\end{align}
and the dual problem consists in maximizing $R$, that is
\begin{align}\label{R}\tag{$\mathcal{P}^*$}
 R\rightarrow\max\text{ in }L^\infty(\Omega,\R^2)\times L^\infty(\Omega,\R^{2\times 2}).
\end{align}
\begin{Satz}\label{Thm1.2}
 Under our general assumptions regarding $\Omega$, $D$, $f$, $F$ and $G$ it holds:
\begin{enumerate}[a)]
 \item Problem $(\mathcal{P}^*)$ has at least one solution $(\rho,\sigma)\in L^\infty(\Omega,\R^2)\times L^\infty(\Omega,\R^{2\times 2})$.
 \item The problems $(\mathcal{P})$ and $(\mathcal{P}^*)$ are related via the so called "inf-sup" relation:
\begin{align*}
 \inf_{\ft{w}\in V}E(\ft{w})=\sup_{\ft{y}\in Y^*}R(\ft{y}),
\end{align*}
i.e. there is no duality gap.
\item Let $(u,v)\in BV(\Omega)\times BV(\Omega,\R^2)$ be a solution of the relaxed problem $(\widetilde{\mathcal{P}})$. Then  the following formula holds:
\begin{align}
 (\rho,\sigma)=D\mathcal{F}\big(\Lambda^a(u,v)\big)=\beta DG(\nabla^au-v)\oplus \alpha DF(\nabla^av)\text{ a.e. on }\Omega,
\end{align}
where we declare 
\[
\Lambda^a(u,v):=(\nabla^a u-v,\nabla^a v). 
\]
In particular, the solution of the dual problem is unique by Theorem \ref{Thm1.1} b).
\end{enumerate}
\end{Satz}

\begin{Rem}
\begin{enumerate}[i)]
\item We would like to advise the reader of the fact, that Theorems \ref{Thm1.1} and \ref{Thm1.2} are valid in arbitrary dimensions $n\geq 2$.
\item  Moreover, both results remain true if we replace the quantity $\intomd (u-f)^2\dx$ with a more general data term $\intomd \Phi\big(|u-f|\big)\dx$ for a strictly convex, increasing and differentiable function 
\[
\Phi:[0,\infty)\rightarrow [0,\infty) 
\]
and consider the problem 
\begin{align}\label{Ephi}\tag{$\widetilde{\mathcal{P}}_\Phi$}
 \begin{split}
\widetilde{E}_\Phi(u,v)=\alpha\intom F(\nabla v)+\beta\intom G(\nabla u-v\cdot\mathcal{L}^2)+\intop_{\Omega-D}\Phi\big(|u-f|\big)\dx\\
  \rightarrow\min\text{ in }BV(\Omega)\times BV(\Omega,\R^2).
\end{split}
\end{align}
\end{enumerate}
\end{Rem}
In order to obtain more regular minimizers, we need to refine our assumptions on the functions $F$ and $G$. In fact, previous work in this regard (see e.g. \cite{BF1}, \cite{BF2}, \cite{BF3} or \cite{FMT} and particularly \cite{BFMT} for more recent results) indicate, that the correct framework for establishing ''classical'' solvability (i.e. in a Sobolev space) of our primal problem \gr{E} is the concept of   ''$\mu$-ellipticity''. This is to say, that we replace the rather general ellipticity condition (F\ref{F2}) with the stronger assumption
\begin{align}
 \label{muell}\tag*{$(\text{F}2)'$}
\begin{split}
 &c_1\frac{1}{(1+|p|)^\mu}|q|^2\leq D^2F(p)(q,q)\leq c_2\frac{1}{1+|p|}|q|^2,\text{ for some }\\
 &c_1,c_2>0,\,\text{ parameter }\mu\in (1,\infty)\text{ and for all }p,q\in\zkz
\end{split}
\end{align}
and (G\ref{G2}) is replaced with
\begin{align}
 \label{munuell}\tag*{$(\text{G}2)'$}
\begin{split}
 &c_1\frac{1}{(1+|y|)^\nu}|x|^2\leq D^2G(y)(x,x)\leq c_2\frac{1}{1+|y|}|x|^2 \text{ for some }\\
 &c_1,c_2>0,\,\text{ parameter }\nu\in (1,\infty)\text{ and for all }x,y\in\R^2.
\end{split}
\end{align} 
We furthermore have to distinguish the case of pure denoising $D=\emptyset$ from the general case. Then we have:
\begin{Satz}
 \label{Thm1.3}
Together with our general assumptions regarding $\Omega$ and $f$, assume $D=\emptyset$ (pure denoising) and  let $F$ satisfy (F1), (F2)', (F3) and let $G$ satisfy (G1), (G2)', (G3) for parameters
\begin{align}\label{mu}
 (\mu,\nu)\in (1,\nicefrac{3}{2})\times (1,2).
\end{align}
Then problem $(\mathcal{P})$ admits a unique solution $(u,v)$ in the Sobolev class 
\[
V=W^{1,1}(\Omega)\times W^{1,1}(\Omega,\R^2). 
\]
It even holds $(u,v)\in W^{1,p}_\loc(\Omega)\times W^{1,p}_\loc(\Omega,\R^2)$ for every $p\in[1,\infty)$.
\end{Satz}
\begin{Rem}\label{Remlin}\begin{enumerate}[i)]
\item The uniqueness of a possible Sobolev-minimizer follows from Theorem \ref{Thm1.1} part c).
\item The results from \cite{BFMT} indicate, that the bound $\mu<\nicefrac{3}{2}$ is not optimal, whereas in \cite{FMT} Fuchs, Tietz and the author have shown, that $\mu,\nu<2$ is indeed necessary for the existence of a solution in the Sobolev class, if $\mu$-elliptic densities are considered.
\item In the case $D\neq \emptyset$, we were not able to prove the above result in full generality. However, if we replace the quadratic error term $\intomd |u-f|^2\dx$ with
\[
 \intom \omega\big(|u-f|\big)dx
 \]
 for a convex, differentiable and increasing function $\omega:\R_0^+\rightarrow\R_0^+$ of linear growth and \gr{mu} with the condition
\begin{align}\tag*{(\ref{mu})'}\label{mu'}
  (\mu,\nu)\in (1,\nicefrac{3}{2})\times (1,\nicefrac{3}{2}),
\end{align}
then the statement of Theorem \ref{Thm1.3} holds for the generalized solution of this modified problem even for $D\neq \emptyset$. Of course, this only concerns the regularity of $u$ and the assumptions of Theorem \ref{Thm1.3} suffice to establish $v\in W^{1,1}(\Omega,\R^2)$  for arbitrary $D$.
\end{enumerate}
\end{Rem}

For regularity results in terms of classical differentiability, we need to put further restraints on our density functions. As our coupled model resembles a vector-valued situation, it is natural to impose a structure condition of Uhlenbeck-type on $F$ in addition to (F1), (F2)' and (F3), which means that we consider functions of the special form
\begin{align}\label{F4}\tag*{(\text{F}4)}
F(p)=g\big(|p|^2\big)
\end{align}
with a convex, increasing function $g:[0,\infty)\rightarrow [0,\infty)$ which is at least of class $C^2$. Again we restrict ourselves to the case of pure denoising ($D=\emptyset$). Then we have:
\begin{Satz}\label{Thm1.4}
Together with our general assumptions regarding $\Omega$ and $f$, assume $D=\emptyset$ and let $F$ satisfy (F1), (F2)', (F3), (F4) and let $G$ satisfy (G1), (G2)' and (G3) with parameter $\mu$ and $\nu$ satisfying \gr{mu}. Let $(u,v)$ be the $W^{1,1}(\Omega)\times W^{1,1}(\Omega,\R^2)$-minimizer from Theorem \ref{Thm1.3}. Then there is an open set $\Omega_0\subset\Omega$ of full measure, i.e.
\begin{align}\label{fm}
 \leb(\Omega-\Omega_0)=0 ,
\end{align}
such that
\[
 (u,v)\in C^{1,\alpha}(\Omega_0)\times C^{1,\beta}(\Omega_0,\R^2) \text{ for all } (\alpha,\beta)\in (0,1)\times (0,1).
\]
For the set $\Omega-\Omega_0$ of possible singularities it further holds 
\begin{align}\label{Hdim}
 \mathcal{H}\text{-}\mathrm{dim}(\Omega-\Omega_0)=0,
\end{align}
which means that the  $\eps$-dimensional Hausdorff measure $\mathcal{H}^\eps(\Omega_0)$ is zero for every $\eps>0$.
\end{Satz}

\begin{Rem}
For $D\neq\emptyset$ the statement of Theorem \ref{Thm1.4} still holds for the modified problem from Remark \ref{Remlin} c) and with \ref{mu'} instead of \gr{mu}. 
\end{Rem}

\begin{Rem}
In contrast to Theorems \ref{Thm1.1} and \ref{Thm1.2}, the statements of Theorems \ref{Thm1.3} and \ref{Thm1.4} crucially depend on the assumption $\Omega\subset\R^2$. 
\end{Rem}
\end{section}

\begin{section}{Relaxation in $BV$, proof of Theorem \ref{Thm1.1}}
As in \cite{FT} or \cite{FM}, a key tool in the proof of Theorem \ref{Thm1.1} is the following density result (cf. Lemma 2.2 in \cite{FT} or Theorem 1.1 in \cite{FM} for a generalization to higher orders):
\begin{Lem}\label{dl}
Let our general assumptions regarding $\Omega$ and $D$ hold. Given $(u,v)\in BV(\Omega)\times BV(\Omega,\R^2)$, there is a sequence $(\varphi_n,\psi_n)\subset C^{\infty}(\overline{\Omega})\times C^\infty(\overline{\Omega},\R^2)$ such that
\begin{gather}
\varphi_n\rightarrow u\text{ in }L^2(\Omega),\label{dl1}\\
\psi_n\rightarrow v\text{ in }L^2(\Omega,\R^2),\label{dl2}\\
\intom\sqrt{1+|\nabla\psi_n|^2}\dx\rightarrow\intom\sqrt{1+|\nabla v|^2},\label{dl3}\\
\intom\sqrt{1+|\nabla\varphi_n-\psi_n|^2}\dx\rightarrow\intom\sqrt{1+|\nabla u-v|^2}\label{dl4}.
\end{gather}
\end{Lem}

\noindent Proof of the Lemma. First we note, that the existence of a sequence $(\psi_n)$ with the properties \gr{dl2} and \gr{dl3} follows directly from Lemma 2.2 in \cite{FT} (note that due to $\Omega\subset\R^2$ it holds $u\in L^2(\Omega)$, $v\in L^2(\Omega,\R^2)$ thanks to embedding theorems). Let us define a linear differential operator with constant coefficients by 
\[
S:C^\infty(\Omega)\times C^\infty(\Omega,\R^2)\rightarrow C^\infty(\Omega,\R^2),\;(\eta,\vartheta)\mapsto \nabla\eta-\vartheta.
\]
The operator $S$ is of local type in the sense of \cite{DT}, p. 688 and we may therefore quote Theorem 2.2 from this work (see also Remark 2.1 therein) to conclude, that there is a sequence $(\varphi_n,\widetilde{\psi}_n)\subset C^\infty(\overline{\Omega})\times C^\infty(\overline{\Omega},\R^2)$ such that
\begin{align}\label{ast}
\begin{split}
&\hspace{4cm}(\varphi_n,\widetilde{\psi}_n)\rightarrow (u,v)\text{ in }L^2(\Omega)\times L^2(\Omega,\R^2)\\
&\text{and }\intom \sqrt{1+|S(\varphi_n,\widetilde{\psi}_n)|^2}\dx\rightarrow \intom  \sqrt{1+|S(u,v)|^2}=\intom  \sqrt{1+|\nabla u-v|^2}.
\end{split}
\end{align}
Furthermore,
\begin{align*} 
\left|\intom \sqrt{1+|S(\varphi_n,\widetilde{\psi}_n)|^2}\dx-\intom \sqrt{1+|S(\varphi_n,\psi_n)|^2}\dx\right|\leq &c\intom |\widetilde{\psi}_n-\psi_n|\dx\rightarrow 0,
\end{align*}
since $\widetilde{\psi}_n,\psi_n\rightarrow v$ in $L^1(\Omega,\R^2)$. Together with \gr{ast}, this proves that $(\varphi_n,\psi_n)$ is a sequence as claimed in the Lemma. \qed

We continue with the proof of Theorem \ref{Thm1.1}. Ad a). Let $(u_k,v_k)\in BV(\Omega)\times BV(\Omega,\R^2)$ denote an $\widetilde{E}$ minimizing sequence. By Lemma \gr{dl} in combination with Reshetnyak's continuity theorem (see, e.g. \cite{AG}, Proposition 2.2) we may assume $(u_k,v_k)\in W^{1,1}(\Omega)\times W^{1,1}(\Omega,\R^2)$. Thanks to the linear growth of $F$ and $G$ it is clear that there are constants $M_1,M_2,M_3>0$ such that
\begin{align}
&\sup_{k\in\N}\intom|\nabla v_k|\dx\leq M_1,\label{m1}\\
&\sup_{k\in\N}\intom |\nabla u_k-v_k|\dx\leq M_2\label{m2}\\
\text{and }&\sup_{k\in\N}\intop_{\Omega-D} |u_k|\dx\leq M_3\label{m3}.
\end{align}
We need the following version of Poincaré's inequality (see \cite{BFW}, Lemma 4.2):
\begin{Lem}\label{Lemm1}
Let $\Omega\subset\R^2$ be a Lipschitz domain and $\rho\in C^1_0(\Omega)$ such that $\intom \rho\dx=1$. Then there is constant $c>0$ which only depends on $\Omega$, such that for any function $u\in W^{1,1}(\Omega)$ it holds
\[
\left\|u-\intom\rho u\dx\right\|_{1;\Omega}\leq c \|\nabla u\|_{1;\Omega}.
\]
\end{Lem}
Now we choose $\rho$ as in the Lemma and such that $\text{spt}(\rho)\subset\Omega-D$ (note that $\Omega-\overline{D}\neq\emptyset$ by our assumption). Then we have
\begin{align}\label{a1}
\sup_{k\in\N}\left|\intom \rho \nabla u_k\dx-\intom \rho v_k\dx\right|\leq \|\rho\|_\infty \sup_{k\in\N}\intom|\nabla u_k-v_k|\dx\leq \|\rho\|_\infty M_2
\end{align}
Furthermore:
\begin{align}\label{a2}
\left|\intom\rho\nabla u_k\dx\right|=\left|\intom \nabla \rho u_k\dx\right|\leq \|\nabla\rho\|_{\infty}\intop_{\Omega-D}|u_k|\dx\leq \|\nabla\rho\|_{\infty} M_3.
\end{align}
Thus, from \gr{a1} and \gr{a2} we conclude
\begin{align*}
\sup_{k\in\N}\left|\intom \rho v_k\dx\right|<\infty
\end{align*}
and \gr{m1} together with Lemma \ref{Lemm1} thus implies
\begin{align*}
\sup_{k\in\N}\|v_k\|_{1,1;\Omega}<\infty.
\end{align*}
But then the boundedness of $v_k$ in $L^1(\Omega,\R^2)$ along with \gr{m2} and \gr{m3} implies (by another application of Poincaré's inequality):
\[
\sup_{k\in\N}\|u_k\|_{1,1;\Omega}<\infty,
\]
so that $(u_k,v_k)$ is bounded in $BV(\Omega)\times BV(\Omega,\R^2)$. By the $BV$-compactness Theorem (see \cite{AFP}, Theorem 3.23 on p. 132), there exists $(u,v)\in BV(\Omega)\times BV(\Omega,\R^2)$ such that
\[
(u_k,v_k)\rightarrow (u,v)\text{ in }L^1(\Omega)\times L^1(\Omega,\R^2) \text{ and }a.e.  
\]
That $(u,v)$ is indeed $\widetilde{E}$-minimal follows immediately since the relaxation $\widetilde{E}$ is lower semicontinuous with respect to $L^1$-convergence by definition (see \cite{AFP}, Remark 5.46 on p. 303).

\noindent The statements of part b) are a mere consequence of the strict convexity of the functions $F$, $G$ and the data fitting term $|u-f|^2$.

\noindent That every $\widetilde{E}$-minimizer is indeed the $L^1$-limit of an $E$-minimizing sequence in the Sobolev class $V$ follows from Lemma \gr{dl} together with $\widetilde{E}_{|V}=E$. That every such limit in $BV$ minimizes $\widetilde{E}$ is a consequence of the above mentioned lower semicontinuity property of the relaxation. It remains to prove that $(\widetilde{\mathcal{P}})$ has a unique solution if $\mathcal{M}\cap V\neq\emptyset$. Assume therefore, that $(u,v)\in V$ minimizes $\widetilde{E}$ and let $(\widetilde{u},\widetilde{v})$  be another element of $\mathcal{M}$. From $\widetilde{E}(u,v)=E(u,v)=\widetilde{E}(\widetilde{u},\widetilde{v})$ and part b) we infer
\begin{align*}
\intom F^\infty\left(\frac{\nabla^s \widetilde{v}}{|\nabla^s\widetilde{v}|}\right)\mathrm{d}|\nabla^s\widetilde{v}|+\intom G^\infty\left(\frac{\nabla^s \widetilde{u}}{|\nabla^s\widetilde{u}|}\right)\mathrm{d}|\nabla^s\widetilde{u}|=0
\end{align*}
and thus $\nabla^s\widetilde{v}=0$ and $\nabla^s\widetilde{u}=0$, which means $(\widetilde{u},\widetilde{v})\in V$. But then b) implies $\nabla\widetilde{v}=\nabla v$ and thereby $v=\widetilde{v}+c$. Further it follows from $\nabla u-v=\nabla\widetilde{u}-\widetilde{v}$ that $\nabla u=\nabla\widetilde{u}+c$, that is $u(x)=\widetilde{u}(x)+c\cdot x+b$, for some $b,c\in\R^2$. Finally $u=\widetilde{u}$ on $\Omega-D$ along with $\mathcal{L}^2(\Omega-D)>0$ requires $b=c=0$, hence $u=\widetilde{u}$ and $v=\widetilde{v}$. \qed
\end{section}

\begin{section}{Duality, proof of Theorem \ref{Thm1.2}}
 Since our arguments follow the ideas in \cite{FT} very closely, the reader will hopefully approve our attempt to give a rather condensed outline of the proof of Theorem \ref{Thm1.2}, referring to the works \cite{FT} or \cite{MT} for the details.

The proof relies on a suitable approximation of  problem $(\mathcal{P})$ through a family of regularizations. To be precise, for $\delta\in (0,1)$ we look at the problem
\begin{align}
\label{Edel}\tag{$\mathcal{P}_\delta$}
\begin{split}
 E_\delta(u,v):=\frac{\delta}{2}\intom |\nabla u|^2+&|\nabla v|^2\dx+E(u,v)\\
&\rightarrow\min\text{ in }W^{1,2}(\Omega)\times W^{1,2}(\Omega,\R^2).
\end{split}
\end{align}
 \begin{Lem}\label{Lem2.1}
  Under our general assumptions regarding  $\Omega$, $f$, $F$ and $G$ it holds:
 \begin{enumerate}[a)]
  \item For any $\delta\in (0,1)$ problem $(\mathcal{P}_\delta)$ admits a unique solution $\ft{u}_\delta=(\udel,\vdel)$ in the space $W^{1,2}(\Omega)\times W^{1,2}(\Omega,\R^2)$.
  \item The family of the $\fudel$'s fulfills:
\begin{align}
 &\sup_{\delta\in (0,1)}\delta\intom |\nabla \udel|^2+|\nabla \vdel|^2\dx<\infty,\label{2.1a}\\
 \text{as well as }&\sup_{\delta\in (0,1)}\intop_{\Omega-D}|\udel|^2\dx<\infty.\label{2.1b}
\end{align}
  \item It holds (not necessarily uniformly with respect to $\delta$!) 
 \begin{align}
\fudel\in W^{2,2}_{\mathrm{loc}}(\Omega)\times W^{2,2}_{\mathrm{loc}}(\Omega,\R^2).
\end{align}
 \end{enumerate}
\end{Lem}
\noindent\textit{Proof of Lemma \ref{Lem2.1}}. \textit{Ad a).} Let $\delta\in (0,1)$ be fixed. Quoting standard results concerning the weak lower semicontinuity of convex functionals on Sobolev spaces, to prove the existence of a minimizer by the direct method, it suffices to show that any $E_\delta$ minimizing sequence is bounded in $W^{1,2}(\Omega)\times W^{1,2}(\Omega,\R^2)$. So let us fix $\delta\in (0,1)$ and denote by $\ft{u}_k=(u_k,v_k)$, $k\in\N$ such a minimizing sequence. By $E_\delta(\ft{u}_k)\leq E_\delta(0,0)=E(0,0)$, it is clear that $|\nabla v_k|$ and $|\nabla u_k|$ are bounded in $L^2(\Omega)$. Furthermore, $f-u_k$ is bounded in $L^2(\Omega-D)$. By Poincaré's inequality we therefore have
\[
\intom \left|u_k(x)-\overline{u_k(x)}\right|^2\dx\leq c\intom |\nabla u_k|^2\dx
\]
where we have set $\overline{u_k(x)}:=\mint_{\Omega-D}u_k(t)\,\mathrm{d}t$. We infer that $u_k$ is bounded in $W^{1,2}(\Omega)$. But then, 
\[
\intom G(\nabla u_k-v_k)\dx\leq E(0,0)
\]
 implies that also $|v_k|$ is bounded in $L^1(\Omega)$ and another application of Poincaré's inequality yields the boundedness of $v_k$ in $W^{1,2}(\Omega,\R^2)$. 

\noindent\textit{Ad b):} this follows immediately from $E_\delta(\ft{u}_\delta)\leq E_\delta(0,0)=E(0,0)$.

\noindent\textit{Ad c):} let $\delta\in (0,1)$ be fixed. The proof of this statement is a standard application of the difference quotient technique to the quadratic variational problems
\begin{align*}
E_\delta(u,v_\delta)\rightarrow\min\text{ in }W^{1,2}(\Omega),\;\text{ with }v_\delta\text{ fixed }
\end{align*}
and
\begin{align*}
E_\delta(u_\delta,v)\rightarrow\min\text{ in }W^{1,2}(\Omega,\R^2),\;\text{ with }u_\delta\text{ fixed,}
\end{align*} respectively. \qed

The core of the proof of Theorem \ref{Thm1.2} now consists in a careful analysis of the convergence behavior of the $\fudel$ as $\delta$ approaches zero. Our claim is, that (at least for a subsequence) $\fudel$ converges in $L^1(\Omega)\times L^1(\Omega,\R^2)$ towards a solution of the relaxed problem $(\widetilde{\mathcal{P}})$, and that 
\begin{align}\label{sigmadel}
\begin{split}
\pmb{\sigma}_\delta=(\rho_\delta,\sigma_\delta):=&\delta\nabla\fudel+D\mathcal{F}(\Lambda \fudel)\\
=&\big[\delta \nabdel+\beta DG(\nabla\udel-\vdel)\big]\oplus \big[ \delta \navdel +\alpha DF(\nabla\vdel)\big]
\end{split}
\end{align}
converges weakly in $L^2(\Omega,\R^2)\times L^2(\Omega,\R^{2\times 2})$ towards a solution of the dual problem $(\mathcal{P}^*)$. 
\begin{Lem}\label{Lem2.2}
The family $\fudel$ is uniformly bounded in the space $V=W^{1,1}(\Omega)\times W^{1,1}(\Omega,\R^2)$. In particular, there is null-sequence $\delta\downarrow 0$ and a function $\overline{\mathbf{u}}=(\overline{u},\overline{v})\in BV(\Omega)\times BV(\Omega,\R^2)$ such that $\fudel\rightarrow\fu$ in $L^1(\Omega)\times L^1(\Omega,\R^2)$ and a.e. as $\delta\downarrow 0$.
\end{Lem}
\noindent\textit{Proof of Lemma \ref{Lem2.2}. }We start with the observation, that due to $E_\delta(\fudel)\leq E_\delta(0,0)=E(0,0)$ and the linear growth of $F$ and $G$  we have the following bounds:
\begin{align}
 &\sup_{\delta\in (0,1)}\intom |\nabla \vdel|\dx\leq M_1',\label{2.5}\\
 &\sup_{\delta\in (0,1)}\intom |\nabla \udel-\vdel|\dx\leq M_2',\text{ as well as }\label{2.6}\\
 &\sup_{\delta\in (0,1)}\intomd |\nabla \udel|^2\dx\leq M_3'\label{2.7}
\end{align}
for constants $M_1',M_2',M_3'>0$. From here on, we may repeat the arguments from the proof of Lemma \ref{Lemm1} to conclude the boundedness of $(u_\delta,v_\delta)$ in $W^{1,1}(\Omega)\times W^{1,1}(\Omega,\R^2)$. The claimed convergence is then seen to be a consequence of the $BV$-compactness theorem. \qed

The next step in the proof of Theorem \ref{Thm1.2} is to show, that the function $\fu$ from Lemma \ref{Lem2.2} in fact minimizes $\widetilde{E}$. Let us fix a null-sequence $\delta\downarrow 0$ as in Lemma \ref{Lem2.2}. By \gr{2.1a}, it holds
\begin{align}\label{2.9}
\delta\nabla \fudel=(\delta\nabdel,\delta\navdel)\rightarrow 0\text{ in }L^2(\Omega,\R^2)\times L^2(\Omega,\zkz)
\end{align}
and since $|D\mathcal{F}|$ is bounded, we have
\begin{align}
 \label{sigbo} \sup_{\delta\in (0,1)}\intom|\pmb{\sigma}_\delta|^2\dx<\infty.
\end{align}
Thus there exists $\pmb{\sigma}\in L^2(\Omega,\R^2)\times L^2(\Omega,\zkz)$ such that 
\begin{align*}
\pmb{\sigma}_\delta\rightharpoondown \pmb{\sigma}\text{ in }L^2(\Omega,\R^2)\times L^2(\Omega,\zkz)\text{ as }\delta\downarrow 0
\end{align*}
at least for another subsequence. Furthermore, setting $\pmb{\tau}_\delta:=D\mathcal{F}(\Lambda\fudel)$, we may assume that there exists $\pmb{\tau}\in L^\infty(\Omega,\R^2)\times L^\infty(\Omega,\zkz)$ such that
\begin{align}
\pmb{\tau}_\delta\overset{*}{\rightharpoondown}\pmb{\tau}.
\end{align} 
 Due to $\pmb{\sigma}_\delta=\delta\nabla\fudel+\pmb{\tau}_\delta$ and \gr{2.9} it must hold
\begin{align}
\pmb{\sigma}=\pmb{\tau}.
\end{align}
Next we observe, that thanks to its minimality with respect to $E_\delta$, $\fudel$ satisfies the following weak Euler-Lagrange equation:
\begin{align}\label{EL}\tag{EL}
\begin{split}
 \delta\intom \nabla\fudel\odot \nabla \pmb{\phi}\dx+ &\intom D\mathcal{F}(\Lambda \fudel)\odot \Lambda\pmb{\phi}\dx+2\intop_{\Omega-D}(u-f)\varphi\dx=0\\
&\text{ for all }\pmb{\phi}=(\varphi,\psi)\in W^{1,2}(\Omega,\R)\times W^{1,2}(\Omega,\R^2).
\end{split}
\end{align}
This can be decoupled into the two equations
\begin{align}
\label{El1}\tag{EL1}
\begin{split}
\intom DF_{\delta}(\nabla v_\delta):\nabla\psi\dx-&\beta\intop_{\Omega}DG(\nabla \udel-\vdel)\cdot\psi\dx=0\\
&\hspace{2cm}\text{ for all }\psi\in W^{1,2}(\Omega,\R^2),
\end{split}
\end{align}
where we have abbreviated $F_{\delta}(p):=\frac{\delta}{2}|p|^2+\alpha F(p)$ for all $p\in \R^{2\times 2}$, and
\begin{align}
\label{El2}\tag{EL2}
\begin{split}
0=\delta \intom \nabla\udel\cdot\nabla \varphi\dx+\beta &\intom DG(\nabla \udel-v_\delta)\cdot\nabla\varphi\dx\\
&+2\intop_{\Omega-D} (\udel-f)\varphi\dx\text{ for all }\varphi\in W^{1,2}(\Omega).
\end{split}
\end{align}
Using the same arguments as in \cite{FT}, Section 4, \gr{EL} along with the duality relation (see \cite{ET}, Proposition 5.1 on p. 21) 
\begin{align*}
\mathcal{F}(\Lambda\fudel)=\pmb{\tau}_\delta\odot \Lambda\fudel-\mathcal{F}^*(\pmb{\tau}_\delta)
\end{align*}
suffices to establish the formula
\begin{align}\label{2.12}
\begin{split}
E_\delta(\fudel)&=-\frac{\delta}{2}\intom |\nabla\fudel|^2\dx-\intom \mathcal{F}^*(\pmb{\tau}_\delta)\dx-\intomd\udel^2\dx+\intomd f^2\dx\\
&\leq -\intom \mathcal{F}^*(\pmb{\tau}_\delta)\dx-\intomd\udel^2\dx+\intomd f^2\dx.
\end{split}
\end{align}
 Note that from the definition of $R$ (see \gr{R1}) it is clear that
\begin{align*}
 \sup_{\mathbf{y}\in V^*}R(\mathbf{y})\leq \inf_{\ft{w}\in V}E(\ft{w})
\end{align*}
and we observe $\inf_{\ft{w}\in V}E(\ft{w})\leq E(\fudel)\leq E_\delta(\fudel)$, so that applying $\limsup\limits_{\delta\downarrow 0}$ on both sides of \gr{2.12} yields
\begin{align}\label{2.13}
  \sup_{\mathbf{y}\in V^*}R(\mathbf{y})\leq \inf_{\ft{w}\in V}E(\ft{w})\leq -\intom \mathcal{F}^*(\pmb{\tau})\dx-\intomd\overline{u}^2\dx+\intomd f^2\dx,
\end{align}
where we have used the convexity of $\mathcal{F}^*$ and Fatou's Lemma.
Following the arguments in \cite{FT}, Section 4 (with $F$ replaced by $\mathcal{F}$ and $\nabla$ replaced by $\Lambda$) we furthermore obtain
\begin{align*}
R(\pmb{\tau})=-\intom\mathcal{F}^*(\pmb{\tau})\dx+\inf_{\mathbf{w}=(w_1,w_2)\in V}\left[\intomd (\overline{u}-w_1)^2+f^2-\overline{u}^2\dx\right]
\end{align*}
and since $\big[\cdots\big]$ is obviously minimal for $w_1=\overline{u}$, we infer from \gr{2.13} the equation
\begin{align*}
 \sup_{\mathbf{y}\in V^*}R(\mathbf{y})=R(\pmb{\tau})=\inf_{\ft{w}\in V}E(\ft{w}),
\end{align*}
i.e. the inf-sup-relation. Further we see that $\pmb{\tau}=\pmb{\sigma}$ maximizes the dual functional and therefore \gr{2.12} implies that (a subsequence of) $\fudel$ is in fact an $E$ minimizing sequence. Part a) and b) of Theorem \ref{Thm1.2} are thus proved.

\noindent For  part c) we claim that it is enough to revise the steps in the proof of Theorem 1.3 in \cite{FT} with $F$ replaced by $\mathcal{F}$ and $\nabla^a$ replaced by $\Lambda^a$, and we therefore would  like to omit the details. \qed
\end{section}

\begin{section}{Sobolev solutions, proof of Theorem \ref{Thm1.3}}
Under the assumptions of Theorem  \ref{Thm1.3}, let $\fudel=(\udel,\vdel)$ be the $\widetilde{E}$-minimizing sequence as constructed in the previous section. Our proof mainly relies on the following lemma:
\begin{Lem}\label{Lem3.1}
Under the assumptions of Theorem \ref{Thm1.3} we have that
\begin{align}
&\varphi_\delta:=\big(1+|\navdel|\big)^{1-\frac{\mu}{2}},\label{1}\\
&\widetilde{\varphi}_\delta:=\big(1+|\nabla \udel|\big)^{1-\frac{\nu}{2}}\label{wthdel}
\end{align}
are uniformly bounded in $W^{1,2}_\loc(\Omega)$. 

\end{Lem}

\noindent\textit{Proof of Lemma \ref{Lem3.1}:} Ad \gr{1}. Throughout, we use summation convention with respect to the index $i\in\{1,2\}$ and denote by $c$ a generic constant. We start with the discussion of the quantity $\phidel$. First we note, that the uniform boundedness of $\varphi_\delta$ in $L^2_\loc(\Omega)$ is clear since we assume $\mu>1$ and $\vdel$ is uniformly bounded in $W^{1,1}(\Omega,\R^2)$ by Lemma \ref{Lem2.2}. Choosing $\palpha\psi$ instead of $\psi$ in the Euler equation  \gr{El1} and performing an integration by parts, we obtain
\begin{align}\label{EL1'}\tag*{(EL1)'}
0=\intom D^2F_\delta(\nabla\vdel)(\palpha\nabla\vdel,\nabla\psi)\dx+\beta\intom DG(\nabu-\vdel)\cdot\palpha\psi\dx,
\end{align}
which, by approximation, holds for all $\psi\in \mathring{W}^{1,2}(\Omega,\R^2)$. % and due to an approximation argument even for $\psi\in W^{1,1}(\Omega,\R^2)$.  
Let $x_0\in\Omega$ be some point and $R>0$ such that $B_{2R}(x_0)\subset\Omega$. We choose $\psi=\eta^2\palpha\vdel$, where $\eta\in C^1_0(\Omega)$ is such that
\begin{align}
\begin{split}\label{eta}
\begin{cases}
&\text{spt}(\eta)\subset B_{2R}(x_0),\\
&0\leq \eta\leq 1\text{ and }\eta\equiv 1\text{ on }B_R(x_0),\\
&|\nabla\eta|\leq\frac{c}{R}.
\end{cases}
\end{split}
\end{align}
\ref{EL1'} then reads
\begin{align*}
0=\intom D^2\Fdel(\nabla\vdel)\big(\palpha\nabla\vdel,\nabla(\eta^2\palpha \vdel)\big)\dx+\beta\intom DG(\nabla\udel-\vdel)\cdot\palpha(\eta^2\palpha\vdel)\dx
\end{align*}
which can be expanded to
\begin{align}\label{3.3}
\begin{split}
0=\intom &D^2\Fdel(\nabla\vdel)\big(\palpha\nabla\vdel,\palpha \nabla\vdel\big)\eta^2\dx\\
&+\underset{\mbox{$=:T_1$}}{\underbrace{2\intom D^2\Fdel(\nabla\vdel)\big(\eta\palpha\nabla\vdel,\nabla\eta\otimes \palpha\vdel\big)\dx}}\\
&\hspace{2cm}+\underset{\mbox{$=:T_2$}}{\underbrace{\beta\intom DG(\nabla\udel-\vdel)\cdot\palpha(\eta^2\palpha\vdel)\dx}}.
\end{split}
\end{align}
We define
\begin{align*}
 \thdel:=D^2\Fdel(\navdel)\big(\palpha\navdel,\palpha\navdel\big)^\frac{1}{2}
\end{align*}
and thus may write \gr{3.3} for short as
\begin{align}\label{3.4}
\intom \thdel^2\eta^2\dx=-T_1-T_2.
\end{align}
Recalling \ref{muell}, we see that the first claim of Lemma \ref{Lem3.1} follows via a uniform estimate of the integral $\intom \thdel^2\eta^2\dx$ on the left-hand side of \gr{3.4}. So let us have a look at the quantity $T_1$ first. Applying the Cauchy-Schwarz inequality to the bilinear form $D^2\Fdel(\nabla\vdel)(\,\cdot\,,\,\cdot\,)$ and then Young's inequality we obtain ($\eps>0$ is arbitrary)
\begin{align}
|T_1|\leq \eps \intom\thdel^2\eta^2\dx+\eps^{-1}\intom D^2\Fdel(\nabla\vdel)\big(\nabla\eta\otimes \palpha\vdel,\nabla\eta\otimes \palpha\vdel\big)\dx.
\end{align}
Choosing $\eps=\frac{1}{2}$, the first summand can be absorbed in the right-hand side of \gr{3.4} whereas to the second summand we apply the estimate \ref{muell} as well as Lemma \ref{Lem2.1} b) and Lemma \ref{Lem2.2} with the result
\begin{align}
\begin{split}
\intom D^2\Fdel(\nabla\vdel)\big(\nabla\eta\otimes \palpha\vdel,\nabla\eta\otimes \palpha\vdel\big)\dx\\
\leq \frac{c}{R^2}+\frac{c}{R^2}\intom\frac{1}{1+|\nabla\vdel|}|\nabla \vdel|^2\dx=c(R).
\end{split}
\end{align}
Hence, we have shown
\begin{align*}
\intom \thdel^2\eta^2\dx\leq c(R)-T_2,
\end{align*}
and it remains to estimate $T_2$. Therefore, we notice that due to our assumption (G1) on the function $G$, we have $DG(\nabu-\vdel)\in L^\infty(\Omega,\R^2)$ uniformly and therefore
\begin{align*}
|T_2|&\leq c\intom |\palpha(\eta^2\palpha\vdel)|\dx\leq c\intom|\nabla\eta||\nabla\vdel|\dx+c\intom \eta^2|\nabla^2\vdel|\dx\\
&\hspace{-0.62cm}\overset{\text{Lemma \ref{Lem2.2}}}{=}c(R)+c\intom \eta^2|\nabla^2\vdel|\dx\\
&=:c(R)+T_3.
\end{align*}
For the quantity $T_3$ we observe
\begin{align*}
 T_3=  \intom \eta^2\big(1+|\nabla \vdel|\big)^\frac{\mu}{2}\frac{|\nabla^2\vdel|}{\big(1+|\nabla \vdel|\big)^\frac{\mu}{2}}\dx,
\end{align*}
which, using Young's inequality can be estimated through 
\begin{align}\label{3.9'}
     T_3\leq \eps\intom\eta^2\frac{|\nabla^2\vdel|^2}{\big(1+|\nabla \vdel|\big)^\mu}\dx+ \eps^{-1}\intom \eta^2\big(1+|\nabla \vdel|\big)^\mu\dx,                                    
\end{align}
where $\eps>0$ is arbitrary.
Therefore \ref{muell} implies
\begin{align*}
 T_3\leq c\eps\intom \thdel^2\eta^2\dx+\eps^{-1}\intom \eta^2\big(1+|\nabla \vdel|\big)^\mu\dx=:c\eps\intom \eta^2\thdel^2\dx+T_4.
\end{align*}
Choosing $\eps$ small enough, we may absorb the first term  in the left-hand side of \gr{3.4}. Further, with 
\begin{align}
\omega_\delta:=\big(1+|\navdel|\big)^{\frac{\mu}{2}}\label{2}
\end{align}
 we may write
\begin{align*}
 T_4=\intom (\eta\psidel)^2\dx.
\end{align*}
Observing the relation $\omega_\delta=\varphi_\delta^{\frac{\mu}{2-\mu}}$, this integral can be treated exactly as the corresponding quantity in \cite{FM}, eq. (4.22) ff., with the result
\begin{align*}
 T_4&\leq c(R)+\intop_{B_{2R}(x_0)} \varphi_\delta^{\frac{4\mu-4}{2-\mu}}\dx\intom \eta^2|\nabla\phidel|^2\dx\\
&\leq c(R)+\intop_{B_{2R}(x_0)} \varphi_\delta^{\frac{4\mu-4}{2-\mu}}\dx\intom \eta^2\thdel^2\dx.
\end{align*}
Here, our assumption (\ref{mu}) is indispensable for $\frac{4\mu-4}{2-\mu}(1-\frac{\mu}{2})=2\mu-2<1$, which enables us to apply Hölder's inequality:
\begin{align}
 \intop_{B_{2R}(x_0)} \varphi_\delta^{\frac{4\mu-4}{2-\mu}}\dx\leq \pi^s R^{2s}\left(\intom 1+|\nabla\vdel|\dx\right)^{2\mu-2}\leq cR^{2s}
\end{align}\label{3.9}
where $s=3-2\mu>0$ and the constant $c$ is independent of the Radius $R$. Combining our estimates of $T_1$ and $T_2$ with \gr{3.4}, we arrive at
\begin{align}
 (1-cR^{2s})\intom \eta^2\thdel^2\dx\leq c(R).
\end{align}
Thus, for radii $R<R_0$ and  $R_0$ such that $cR_0^{2s}<1$, we have the uniform estimate
\begin{align}\label{3.12'}
 \intom \eta^2\thdel^2\dx\leq c(R_0).
\end{align}
Claim \gr{1} of Lemma \ref{Lem3.1} now follows from a covering argument. 

Note, that as a consequence of  \gr{1} and Sobolev's embedding Theorem (recall $n=2$), we have
\begin{align}\label{Lp}
\nabla \vdel\in L^p_\mathrm{loc}(\Omega,\zkz)\text{ for any }p\in[1,\infty),
\text{ uniformly with respect to $\delta$.}
\end{align}
In particular, 
\begin{align}
\label{Linfty} \vdel \in L_\loc^\infty(\Omega,\zkz)\text{ uniformly with respect to $\delta$.}
\end{align}
Furthermore, it follows (at least after passing to a suitable subsequence $\delta\downarrow 0$) that $\nabla\vdel$ has a weak $L^p_\loc(\Omega,\R^2)$-limit for some $p>1$ and since $\vdel\rightarrow v$ in $L^1(\Omega,\R^2)$ and a.e., we infer $v\in W^{1,p}_\loc(\Omega,\R^2)$. Eventually, $v\in W^{1,1}_\mathrm{loc}(\Omega,\R^2)\cap BV(\Omega,\R^2)$ which is a subset of  $W^{1,1}(\Omega,\R^2)$.

Let us now turn to the corresponding quantity $\widetilde{\varphi}_\delta$ involving $\udel$.
We start with the Euler equation \gr{El2} (keep in mind that $D=\emptyset$ in the setting of Theorem \ref{Thm1.3}), where we choose $\varphi=\palpha(\eta^2\palpha\udel)$ for some $\eta\in C_0^1(\Omega_0)$ satisfying the set of conditions from \gr{eta}. Writing $G_\delta(x):=\delta  |x|^2+\beta G(x)$ for $x\in\R^2$, \gr{El2} reads after an integration by parts:
\begin{align}\label{above}
\begin{split}
\intom D^2G_\delta&(\nabla\udel-\vdel)\big(\palpha\nabla\udel,\nabla(\eta^2\palpha\udel)\big)\dx\\
&-\beta\intom D^2G(\nabla\udel-\vdel)\big(\palpha\vdel,\nabla(\eta^2\palpha\udel)\big)\dx\\
&\hspace{2cm}-2\intom(\udel-f)\palpha(\eta^2\palpha\udel)\dx=0.
\end{split}
\end{align}
Now we  define
\[
 \widetilde{\Theta}_\delta:=D^2G_\delta(\nabla\udel-\vdel)\big(\palpha\nabla\udel,\pnu\big)^{1/2}
\]
due to which we may write \gr{above} as
\begin{align}\label{3.20}
\begin{split}
 \intom \eta^2\wthdel^2\dx=&-2\intom D^2G_\delta(\nabla\udel-\vdel)\big(\eta\palpha\nabla\udel,\nabla\eta\palpha\udel\big)\dx\\
  &+\beta\DG\big(\eta\palpha\vdel,\eta\palpha\nabla\udel\big)\dx\\
  &+2\beta\DG\big(\eta\palpha\vdel,\nabla\eta\palpha\udel\big)\dx\\
  &+2\intop_{\Omega}(\udel-f)\palpha(\eta^2\palpha\udel)\dx\\
  &=:T_1'+T_2'+T_3'+T_4'.
\end{split}
\end{align}
First, we note that due to \ref{munuell}  it holds
\begin{align}\label{wthph}
\begin{split}
\eta^2 \wthdel^2&\geq c_1 \frac{1}{\big(1+|\nabla\udel-\vdel|\big)^{\nu}}|\eta\nabla^2\udel|^2\\
 &\geq c_1 \frac{1}{\big(1+|\nabla\udel|+|\vdel|\big)^\nu}|\eta\nabla^2\udel|^2\\
 &\overset{\gr{Linfty}}{\geq} \widetilde{c}\frac{1}{\big(1+|\nabla\udel|\big)^\nu}|\eta\nabla^2\udel|^2\\
 &\geq c_2 \eta^2|\nabla\phitil|^2.
 \end{split}
\end{align}
Hence, by \gr{3.20} and our choice of $\eta$ we have
\begin{align*}
 \intop_{B_R(x_0)} |\nabla\phitil|^2\dx\leq c\big(T_1'+T_2'+T_3'+T_4'\big).
\end{align*}
In the Integral $T_1'$, we first apply the Cauchy-Schwarz inequality to the bilinear form $D^2G_\delta(\nabla\udel-\vdel)(\cdot,\cdot)$, followed by Young's inequality to obtain
\begin{align*}
 |T_1'|\leq \eps &\intom\eta^2\wthdel^2\dx+c(\eps)\intom D^2G(\nabu-\vdel)(\nabla\eta\palpha\udel,\nabla\eta\palpha\udel)\dx,
\end{align*}
where $\eps>0$ is arbitrarily small.
The first summand can be absorbed in the left-hand side of \gr{3.20}. For the second term, we consider the set
\begin{align*}
 \Sigma:=\left\{x\in B_{2R}(x_0)\,:\,|\nabla\udel|\leq|\vdel|+1\right\}.
\end{align*}
We observe that
\begin{align*}
 &\DG(\nabla\eta\palpha\udel,\nabla\eta\palpha\udel)\dx\\
 &\hspace{1cm}\leq\intop_{\Sigma}|D^2G||\nabla\eta|^2|\nabla\udel|^2\dx+\intop_{B_{2R}(x_0)-\Sigma}D^2G(\nabla\udel-\vdel)(\nabla\eta\palpha\udel,\nabla\eta\palpha\udel)\dx\\
 &\hspace{0.6cm}\overset{\text{\ref{munuell}}\,\&\,\gr{Linfty}}{\leq}c\left( \frac{1}{R^2}+\frac{1}{R^2}\intop_{B_{2R}(x_0)-\Sigma}\frac{1}{1+|\nabla\udel-\vdel|}|\nabla\udel|^2\dx\right)\\
 &\hspace{1cm}\leq c\left( \frac{1}{R^2}+\frac{1}{R^2}\intop_{B_{2R}(x_0)-\Sigma}\frac{1}{1+|\nabla\udel|-|\vdel|}|\nabla\udel|^2\dx\right)\leq c(R).
 \end{align*}
To the quantity $T_2'$ we apply the Cauchy-Schwarz inequality and then Young's inequality with the following result:
\begin{align*}
 |T_2'|\leq \eps\intom\eta^2\wthdel^2\dx+c(\eps)\DG(\eta\palpha\vdel,\eta\palpha\vdel)\dx.
\end{align*}
Again, we absorb the first term in the left-hand side of \gr{3.20} and the second term is bounded by (\ref{Lp}). Combining the arguments for $T_1'$ and $T_2'$, we can estimate $T_3'$ by
\begin{align*}
 |T_3'|\leq \frac{c}{R^2}
\end{align*}
and \gr{3.20} reads
\begin{align}\label{3.21}
 \intom \eta^2\wthdel^2\dx\leq c\left(1+\frac{1}{R^2}\right)+|T_4'|.
\end{align}
It thus remains to give a bound for $T_4'$. An integration by parts yields
\begin{align*}
T_4'= -\intom \eta^2 |\nabla\udel|^2\dx-\intom f\palpha(\eta^2\pnu)\dx.
\end{align*}
The Dirichlet integral can be moved to the left-hand side of \gr{3.21}, so that
\begin{align}\label{3.22}
 \intom \eta^2\wthdel^2\dx+\intom\eta^2|\nabla\udel|^2\dx\leq c(R)+\intom |f||\palpha(\eta^2\pnu)|\dx.
\end{align}
Note that by our assumptions we have $f\in L^\infty(\Omega)$ and thus \gr{3.22} together with Lemma \ref{Lem2.2} implies
\begin{align*}
 \intom \eta^2\wthdel^2\dx+\intom\eta^2|\nabla\udel|^2\dx\leq c(R)+c\intom \eta^2|\nabla^2\udel|\dx.
\end{align*}
The non-constant term on the right-hand side can now be estimated just like the corresponding term $T_3$ in  \gr{3.9'}, which yields
\begin{align*}
 \intom \eta^2|\nabla^2\udel|\dx\leq c\eps\intom \eta^2\wthdel^2\dx+\eps^{-1}\intom \eta^2\big(1+|\nabdel|\big)^\mu\dx.
\end{align*}
For $\eps$ small enough, the first term can be absorbed in the left-hand side of \gr{3.22} and to the second term we apply Young's inequality once again (making use of $\mu<2$) which results in
\[
 \intom \eta^2\big(1+|\nabdel|\big)^\mu\dx\leq \eps\intom \eta^2|\nabla u|^2\dx+c(R),
\]
and the Dirichlet integral can be absorbed in the left-hand side of \gr{3.22} (provided $\eps$ is chosen small enough).
Then claim \gr{wthdel} follows from \gr{3.21} and \gr{wthph}.  Via Sobolev's embedding Theorem, \gr{wthdel} yields
\begin{align}\label{Lpu}
 \nabla \udel\in L^p_\mathrm{loc}(\Omega)\text{ for any }p\in[1,\infty)\text{ and uniform with respect to $\delta$},
\end{align}
which allows us to infer $u\in W^{1,1}(\Omega)$ and $u\in W^{1,p}_\loc(\Omega)$. It even follows from \gr{Lpu}: 
\begin{align}
\label{Luinfty}\udel\in L^\infty_\mathrm{loc}(\Omega)\text{ uniform with respect to $\delta$}.
\end{align}

\begin{Rem}
 If $D\neq\emptyset$, we cannot readily perform an integration by parts to estimate the crucial quantity
\[
 T_4':=\intop_{\Omega-D}(\udel-f)\palpha(\eta^2\palpha\udel)\dx.
\]
However, switching to an error term of linear growth as proposed in Remark \ref{Remlin} turns $T_4'$ into
\[
T_4'':=\intop_{\Omega-D}\omega'\big(|\udel-f|\big)\frac{\udel-f}{|\udel-f|}\palpha(\eta^2\palpha\udel)\dx
\]
and since $|\omega'|$ is bounded we may estimate
\begin{align*}
 |T_4''|\leq c\intop_{\Omega-D}|\palpha(\eta^2\palpha\udel)|\dx\leq c(R)+c\intom \eta^2|\nabla^2\udel|\dx.
\end{align*}
Now, employing the same arguments that have been used for the term $T_3$,  we can establish $u\in W^{1,1}(\Omega)$ even for $D\neq\emptyset$. 
\end{Rem}

\end{section}

\begin{section}{Hölder continuity, proof of Theorem \ref{Thm1.4}}
Our proof follows the ideas in \cite{BF2}, which are based on results by Frehse and Seregin from the works \cite{Fr} and \cite{FS}. An essential condition for the application of these techniques is the validity of following lemma: 

\begin{Lem}
 \label{Lem4.1}
Under the assumptions of Lemma \ref{Lem3.1} it holds for any $s\in (1,2)$:
\begin{align}
 \vdel\text{ is uniformly bounded in } W_\loc^{2,s}(\Omega,\R^2)\label{4}
\end{align}
and
\begin{align}
 \udel\text{ is uniformly bounded in } W_\loc^{2,s}(\Omega).\label{5}
\end{align}
Moreover, 
\begin{align}
\omega_\delta:=\big(1+|\navdel|\big)^{\frac{\mu}{2}}\text{ is uniformly bounded in } W_\loc^{1,2}(\Omega)\label{6}
\end{align}
and 
\begin{align}
\widetilde{\omega}_\delta:=\big(1+|\nabdel|\big)^{\frac{\nu}{2}}\text{ is uniformly bounded in } W_\loc^{1,2}(\Omega).\label{7}
\end{align}
\end{Lem}

\noindent Ad \gr{4}. Recalling the definition of $\varphi_\delta$ from Lemma \ref{Lem3.1} as well as inequality \ref{muell}, we see that the uniform boundedness of $\nabla\varphi_\delta$ in $L_\loc^2(\Omega,\R^2)$, which is obtained from \gr{Lp},  implies that for any compact subset $\Omega^*\Subset\Omega$, there is a constant $c(\Omega^*)>0$ (independent of $\delta$!) such that
\begin{align*}
 \intop_{\Omega^*}\frac{|\nabla^2\vdel|^2}{\big(1+|\nabla\vdel|\big)^\mu}\dx\leq c(\Omega^*).
\end{align*}
Let now $s\in(1,2)$ be arbitrary. We may write
\begin{align*}
\intop_{\Omega^*}|\nabla^2\vdel|^s\dx=\intop_{\Omega^*}\left(\frac{|\nabla^2\vdel|^2}{\big(1+|\nabla\vdel|\big)^\mu}\right)^\frac{s}{2}\big(1+|\navdel|\big)^{\mu\frac{s}{2}}\dx
\end{align*}
and an application of Hölder's inequality yields
\begin{align*}
\intop_{\Omega^*}|\nabla^2\vdel|^s\dx\leq \left(\intop_{\Omega^*}\frac{|\nabla^2\vdel|^2}{\big(1+|\nabla\vdel|\big)^\mu}\dx\right)^\frac{s}{2}\left(\intop_{\Omega^*}\big(1+|\navdel|\big)^\frac{2-s}{2}\dx\right)^\frac{2}{2-s},
\end{align*}
so that \gr{4} follows from \gr{1} and \gr{Lp}. The same argument works for \gr{5}.

We continue with \gr{6}. Choose $B_{2R}(x_0)\subset \Omega$ and $\eta$  according to \gr{eta}.

 Setting $\Gamma_\delta:=1+|\navdel|^2$, we observe
\begin{align*}
\intb |\nabla\omega_\delta|^2\dx&=\intb \big(1+|\navdel|\big)^{\mu-2}|\nabla^2\vdel|^2\dx\\
&=\intb\big(1+|\navdel|\big)^{-\mu}|\nabla^2\vdel|^2\big(1+|\navdel|\big)^{2\mu-2}\dx\\
&\overset{\mbox{\ref{muell}}}{\leq}c\intdb\eta	^2D^2F_\delta(\navdel)\big(\palpha \navdel,\palpha\navdel\big)\Gamma_\delta^{\mu-1}\dx.
\end{align*}
Choosing $\psi=\palpha(\eta^2\palpha\vdel\Gamma_\delta^{\mu-1})$ in \gr{El1} yields
\begin{align*}
\intdb\eta^2	&D^2F_\delta(\navdel)\big(\palpha \navdel,\palpha\navdel\big)\Gamma_\delta^{\mu-1}\dx\\
= -&\beta\intdb DG(\nabdel-\vdel)\cdot \palpha(\eta^2\palpha\vdel\Gamma_\delta^{\mu-1})\dx\\
&-\intdb D^2F_\delta(\navdel)\big(\palpha \navdel,\palpha\vdel\otimes\nabla\eta^2\big)\Gamma_\delta^{\mu-1}\dx\\
&-\intdb D^2F_\delta (\navdel)\big(\palpha\navdel,\palpha\vdel\otimes \nabla \Gamma_\delta^{\mu-1}\big)\eta^2\dx\\
&=:-I_1-I_2-I_3.
\end{align*}
We start with the term $I_2$. Applying the Cauchy-Schwarz inequality to the bilinear form $D^2F_\delta (\cdot,\cdot)$ and then Young's inequality yields
\begin{align*}
|I_2|\leq \frac{1}{2}\intdb D^2F_\delta(\navdel)&(\palpha\navdel,\palpha\navdel)\eta^2\dx\\
&+\frac{1}{2}\intdb D^2F_\delta(\navdel)(\palpha\vdel\otimes\nabla\eta,\palpha\vdel\otimes\nabla\eta)\Gamma_\delta^{2\mu-2}\dx\\
&\hspace{-4cm}\leq\frac{c}{R^2}\left[ \intdb \thdel^2\eta^2\dx+\intdb\delta |\navdel|^2\Gamma_\delta^{2\mu-2}+\frac{1}{1+|\navdel|}|\navdel|^2\Gamma_\delta^{2\mu-2}\dx\right]
\end{align*}
and this term is bounded due to \gr{3.12'} and \gr{Lp}. Let us continue with $I_3$. At this point, we make use of the structure condition \ref{F4} which enables us to write (cf. the calculation on the bottom of page 62 in \cite{Bi})
\begin{align*}
D^2F_\delta (\navdel)&\big(\palpha\navdel,\palpha\vdel\otimes \nabla \Gamma_\delta^{\mu-1}\big)\eta^2\\
&=\frac{1}{2}D^2F_\delta (\navdel)\big(e_i\otimes \nabla|\nabla\vdel|^2,e_i\otimes \nabla|\nabla\vdel|^2\big)\Gamma_\delta^{\mu-2}>0,
\end{align*}
where $e_i$ denotes the canonical basis of $\R^2$. Hence we may just neglect the term $I_3$ and it remains to give a bound on the quantity $I_1$. We note, that due to the boundedness of $DG(\nabdel-\vdel)$ it holds
\begin{align}\label{4.14}
\begin{split}
&|I_1|\leq c\intdb  |\palpha(\eta^2\palpha\vdel\Gamma_\delta^{\mu-1})|\dx\\
&\hspace{-0.3cm}\leq \frac{c}{R}\left[\intdb |\navdel|\Gamma_\delta^{\mu-1}\dx+\intdb \eta^2|\nabla^2\vdel|\Gamma_\delta^{\mu-1}\dx\right].
\end{split}
\end{align}
The first term in the bracket is bounded by \gr{Lp}. For the second one, we note that an application of Young's inequality yields
\begin{align*}
 \intdb \eta^2|\nabla^2\vdel|\Gamma_\delta\dx&\leq c\left[\intdb\eta^2 \frac{|\nabla^2\vdel|^2}{\big(1+|\vdel|\big)^\mu}\dx+\intdb\eta^2\Gamma_\delta^{\frac{5}{2}\mu-2}\dx\right]\\
 &\overset{\mbox{\ref{muell}}}{\leq}  c\left[\intom \thdel^2\eta^2\dx+\intdb\Gamma_\delta^{\frac{5}{2}\mu-2}\dx\right]
\end{align*}
and this is bounded due to \gr{3.12'} and \gr{Lp}. Thus, \gr{6} follows. For \gr{7}, we only  note that this follows from similar  arguments and thereby finish the proof of Lemma \ref{Lem4.1}. 

We continue with the proof of Theorem \ref{Thm1.4}. In the differentiated Euler-Lagrange equation \ref{EL1'}, we  now consider $\psi=\eta^2\big(\palpha\vdel-\overline{\palpha\vdel}\big)$, where we set $\overline{\palpha\vdel}:=\mint_\Omega\palpha\vdel\dx$ and $\eta\in C_0^1(\Omega)$ is chosen according to \gr{eta}. Denoting by $T$ the annulus $B_{2R}(x_0)-B_R(x_0)$ (remember $\eta\equiv\text{const.}$ outside $T$), \ref{EL1'} reads
\begin{align*}
 0=\intom\thdel^2 \eta^2\dx+2\intop_{T} &D^2\Fdel(\nabla\vdel)\big(\palpha\nabla\vdel,\nabla\eta\otimes(\palpha\vdel-\overline{\palpha\vdel})\big)\eta \dx\\
 &+\beta \intom \underset{\mbox{$\in L^\infty(\Omega,\R^2)$}}{\underbrace{DG(\nabla\udel-\vdel)}}\cdot\palpha\big(\eta^2(\palpha\vdel-\pvd)\big)\dx
\end{align*}
and we infer, that for some constant $c>0$ independent of $\eta$ it holds
\begin{align}\label{3.11}
\begin{split}
 &\intom \thdel^2 \eta^2\dx\leq c\left[\intop_{T}\big|D^2\Fdel(\nabla\vdel)\big(\palpha\nabla\vdel,\nabla\eta\otimes(\palpha\vdel-\overline{\palpha\vdel})\big)\eta\big| \dx\right.\\
  &\hspace{3.5cm}\left.+\beta\intom\Big|\palpha\big(\eta^2(\palpha\vdel-\pvd)\big)\Big|\dx\right]=:c\bigg[S_1+S_2\bigg].
\end{split}
\end{align}
 In $S_1$ we apply the Cauchy-Schwarz inequality to the bilinear form $D^2\Fdel(\cdot,\cdot)$ and obtain
\begin{align*}
&S_1\leq\intop_{T}\overset{\mbox{$=\thdel$}}{\overbrace{D^2F_\delta(\vdel)\big(\palpha\navdel,\palpha\navdel\big)^\frac{1}{2}}}\\
&\hspace{1.5cm}\cdot D^2F_\delta(\navdel)\big(\nabla\eta\otimes(\palpha\vdel-\overline{\palpha\vdel}),\nabla\eta\otimes(\palpha\vdel-\overline{\palpha\vdel})\big)^\frac{1}{2}\dx.
\end{align*}
Applying Hölder's inequality yields
\begin{align*}
S_1&\leq \left[\intop_{T}\thdel^2\dx\right]^\frac{1}{2}\left[\intop_{T}D^2F_\delta(\navdel)\big(\nabla\eta\otimes(\palpha\vdel-\overline{\palpha\vdel}),\nabla\eta\otimes(\palpha\vdel-\overline{\palpha\vdel})\big)\dx\right]^\frac{1}{2}\\
&\leq c\left[\intop_{T}\thdel^2\dx\right]^\frac{1}{2}\left[\intop_{T}|D^2F_\delta(\navdel)||\nabla\eta|^2|\palpha\vdel-\overline{\palpha\vdel}|^2\dx\right]^\frac{1}{2}
\end{align*}
and since $|D^2F_\delta|$ is bounded, we arrive at
\begin{align}\label{3.12}
\begin{split}
S_1&\leq \frac{c}{R}\left[\intop_{T}\thdel^2\dx\right]^\frac{1}{2}\left[\intop_{T}|\palpha\vdel-\overline{\palpha\vdel}|^2\dx\right]^\frac{1}{2}\\
&\leq \frac{c}{R}\left[\intop_{T}\thdel^2\dx\right]^\frac{1}{2}\intop_T|\nabla^2\vdel|\dx,
\end{split}
\end{align}
where the Sobolev-Poincaré inequality has been applied in the last step.
Next, we note that  by \ref{muell} we have
\begin{align*}
\big(1+|\navdel|\big)^{-\frac{\mu}{2}}|\nabla^2\vdel|\leq c\thdel
\end{align*}
and thus
\begin{align}\label{3.13'}
|\nabla^2\vdel|\leq c\thdel \big(1+|\navdel|\big)^\frac{\mu}{2}=c\thdel\omega_\delta
\end{align}
(with $\omega_\delta$ as in Lemma \ref{Lem4.1}). Consequently,  it follows from \gr{3.12} that
\begin{align}
S_1\leq \frac{c}{R}\left[\intop_T\thdel^2\dx\right]^{1/2}\intop_T|\nabla^2\vdel|\dx\overset{\gr{3.13'}}{\leq}  \frac{c}{R}\left[\intop_T\thdel^2\dx\right]^{1/2}\intop_T\thdel\psidel\dx.
\end{align}
The term $S_2$ can be treated in the same way as the corresponding quantity in \cite{BF2} on p. 164 with the result
\begin{align*}
S_2\leq \intop_T\thdel\psidel\dx+\intom \eta^2\thdel\psidel\dx.
\end{align*}
The estimates of $S_1$ and $S_2$ together with \gr{3.11} now establish the crucial inequality (3.17) from \cite{BF2} in our setting:
\begin{align}\label{crucial}
 \intdb \thdel^2\dx\leq \frac{c}{R}\left[\intop_{T}\thdel^2\dx +R^2\right]^\frac{1}{2}\intop_{T}\thdel\psidel\dx+c\intdb\eta^2\thdel\psidel\dx,
\end{align}
which holds for all radii $0<R<R_0$ and all points $x_0\in\Omega$ such that $B_{2R_0}(x_0)\Subset\Omega$ with a constant $c$ only depending on $R_0$. To the last term, we apply Young's and Hölder's inequality to get
\begin{align*}
\intdb\eta^2\thdel\psidel\dx\leq \frac{1}{2}\intdb\thdel^2\eta^2\dx +cR^{2\frac{q-1}{q}},
\end{align*}
where due to \gr{Lp}, the exponent $q$ can be chosen from $(1,\infty)$. We fix $\gamma<1$ and thus obtain:
\begin{align}\label{crucial2}
 \intb \thdel^2\dx\leq \frac{c}{R}\left[\intop_{T}\thdel^2\dx +R^2\right]^\frac{1}{2}\intop_{T}\thdel\psidel\dx+c R^\gamma.
\end{align}
As it is explained in detail in \cite{BF2}, pp. 164 ff., this inequality suffices to deduce the following growth estimate for the quantity $\thdel$:
\begin{align}\label{3.17}
\intop_{B_R(x_0)}\thdel^2\dx\leq c\frac{1}{\ln\left(\frac{1}{R}\right)^t}\text{ for all }t\geq 1,
\end{align}
 for all balls $B_R(x_0)$ as above with $0<R<R_0$ and with a local constant $c$  only depending on $R_0$.
 Observing, that for $\sigma_\delta=D\Fdel(\navdel)$ (see \gr{sigmadel}) it holds
\begin{align}
\begin{split}\label{sig}
|\nabla\sigma_\delta|^2&=\palpha\sigma_\delta:\palpha\sigma_\delta=D^2\Fdel(\navdel)(\palpha\navdel,\palpha\sigma_\delta)\\
&\leq \bigg(D^2\Fdel(\navdel)(\palpha\navdel,\palpha\navdel)\bigg)^{1/2}\bigg(D^2\Fdel(\navdel)(\palpha\sigma_\delta,\palpha\sigma_\delta)\bigg)^{1/2}\\
&\leq c\thdel|\nabla\sigma_\delta|
\end{split}
\end{align}
and hence 
\[
|\nabla\sigma_\delta|\leq c\thdel,
\]
the estimate \gr{3.17} implies
\[
\intop_{B_R(x_0)}|\nabla\sigma_\delta|^2\dx\leq c\frac{1}{\ln\left(\frac{1}{R}\right)^t}\text{ for all }t\geq 1.
\]
 Along with Lemma \ref{Lem3.1} and \gr{6}, this is enough to infer the continuity of the $\sigma_\delta$ on every ball $B_R(x_0)$ with $R<R_0$ from the results in \cite{Fr}, p. 287 (see also Lemma 6 and 7 in \cite{BF4}), with modulus of continuity given by
\begin{align}\label{modofco}
\sup_{x,y\in B_R(x_0)}|\sigma_\delta(x)-\sigma_\delta(y)|\leq K|\ln R|^{1-\frac{t}{2}},
\end{align}
with a constant $K=K(R_0)$. The uniform boundedness of $\sigma_\delta$ in $L^2(\Omega,\zkz)$ (cf. \gr{sigbo}) along with \gr{modofco} now implies
\begin{align}\label{Linf}
\sup_{\delta\in (0,1)}\|\sigma_\delta\|_{L^\infty(B_R(x_0))}<\infty.
\end{align}
Furthermore, \gr{modofco} yields the equicontinuity of the $\sigma_\delta$ on any compact subset $\Omega^*\Subset\Omega$, such that $B_{2R_0}(x)\subset\Omega$ for all $x\in\Omega^*$. An application of the Arzelà-Ascoli compactness-theorem thus gives  the existence of a continuous function $\sigma$ such that
\begin{align*}
\sigma_\delta\rightarrow\sigma\text{ locally uniformly, }
\end{align*}
at least for a subsequence $\delta\downarrow 0$. By \gr{4}, we may in addition assume $\navdel\rightarrow \nabla v $ a.e., where $(u,v)\in W^{1,1}(\Omega)\times W^{1,1}(\Omega,\R^2)$ is the unique $E$-minimizer. Hence it holds
\begin{align}\label{3.18}
 DF(\nabla v(x))=\sigma(x)
\end{align}
for almost all $x\in\Omega$. We note, that by the inverse function theorem and \ref{muell}, $\text{Im}(DF)$ is an open set. In particular, $\sigma^{-1}(\text{Im}(DF))$ is open and for every point $x_0\in\Omega$, for which \gr{3.18} holds, there is a small ball $B_\eps(x_0)\subset \R^2$ such that $\sigma(x)\in \text{Im}(DF)$ for all $x\in B_\eps(x_0)$. Hence $DF^{-1}(\sigma)$ is a continuous representative of $\nabla v$ on $B_\eps(x_0)$. But due to the continuity of $\sigma$ (and the Lipschitz continuity of $DF$), \gr{3.18} particularly holds for all Lebesgue points  of $\nabla v$. Identifying $\nabla v$ with its Lebesgue point representative, we thus obtain from \gr{3.18} by inversion a continuous representative of $\nabla v$ on the set:
\[
 \Omega_0:=\left\{x\in\Omega\,:\,\lim\limits_{r\downarrow 0}\mint_{B_r(x)}\nabla v\dx\text{ exists in }\R^{2\times 2}\right\},
\]
which alongside is proved to be open. In particular, $\nabla v\in L^\infty_{\loc}(\Omega_0,\zkz)$ and we can therefore argue just like  in \cite{BFW} on p. 76, to  deduce the Hölder continuity of $v$ on $\Omega_0$ from  the hole-filling technique applied to the inequality \gr{crucial2}. That $\Omega-\Omega_0$ does indeed have Hausdorff-dimension $0$ is a consequence of Theorem 2.1 on p. 100 of \cite{Gia} and $v\in W^{2,s}_\loc(\Omega,\R^2)$, $s\in [1,2)$.

We now come to the corresponding statements concerning $u$. In the following calculations, we restrict ourselves to the open subset $\Omega_0\subset\Omega$ on which we have already established local Hölder-continuity of $v$. We introduce a new sequence $(\widetilde{\udel})$ of $\delta$-regularizers which solve
\begin{align*}
\beta \intomo G(\nabla w-v)\dx+\intomod (w-f)^2\dx+\frac{\delta}{2}\intomo |\nabla w|^2\dx\rightarrow\min\text{ in }W^{1,2}(\Omega),
\end{align*}
where $v$ is the Hölder continuous minimizer from above. We note, that due to
\begin{align*}
&E(u,v)\leq \alpha\intom F(\nabla v)\dx+\beta\intom G(\nabla \wudel-v)\dx+\intom (\wudel-f)^2\dx\\
&\leq \alpha\intom F(\nabla v)\dx+\beta\intom G(\nabla \wudel-v)\dx+\intom (\wudel-f)^2\dx+\frac{\delta}{2}\intom |\nabla \wudel|^2\dx\\
&\leq \alpha\intom F(\nabla v)\dx+\beta\intom G(\nabla \udel-v)\dx+\intom (\udel-f)^2\dx+\frac{\delta}{2}\intom |\nabla \udel|^2\dx\\
&\quad \overset{\delta\downarrow 0}{\longrightarrow}E(u,v),
\end{align*}
the sequence $(\wudel,v)$ is $E$-minimizing and Theorem \ref{Thm1.1} c) implies
\[
\wudel\rightarrow u\text{ in }L^1(\Omega)\text{ and a.e. (at least for a subsequence $\delta\downarrow 0$).}
\]
Moreover, we can verify the properties from Lemma \ref{Lem3.1} and Lemma \ref{Lem4.1} for the sequence $(\wudel)$:
\begin{Lem}\label{Lem4.2}
It holds uniformly with respect to the parameter $\delta$:
\begin{align*}
&\wudel\in W^{2,s}_\loc(\Omega_0)\text{ for all }s\in (1,2),\\
&\widehat{\varphi}_\delta:=\big(1+|\nabla\wudel|\big)^{1-\frac{\nu}{2}}\in W^{1,2}_\loc(\Omega_0),\\
&\widehat{\Theta}_\delta:=D^2G_\delta(\nabla\wudel-v)\big(\palpha\nabla\wudel,\palpha\nabla\wudel\big)^{\frac{1}{2}}\in W^{1,2}_\loc(\Omega_0),\\
&\widehat{\omega}_\delta:=\big(1+|\nabla\wudel|\big)^\frac{\nu}{2}\in W^{1,2}_\loc(\Omega_0).
\end{align*}
\end{Lem}
These statements can be proved just like the corresponding results from Lemma \ref{Lem3.1} and Lemma \ref{Lem4.1}, and we do not want to repeat the technical details here.

Continuing with the proof of Theorem \ref{Thm1.4}, we find that $\wudel$ satisfies the following Euler equation:
\begin{align*}
\delta\intomo \palpha\nabla\wudel\cdot\nabla \varphi\dx +&\beta \intomo D^2G(\nabla\wudel-v)\big(\palpha\nabla\wudel-\palpha v,\nabla \varphi\big)\dx\\
&\hspace{2cm}-\intomod(\wudel-f)^2\palpha\varphi=0\text{ for all }\varphi\in \mathring{W}^{1,2}(\Omega_0).
\end{align*}
With the choice $\varphi:=\eta^2(\palpha\wudel-\overline{\palpha\wudel})$ (now with $\eta\in C^\infty_0(\Omega_0)$ and the properties \gr{eta}), we deduce
\begin{align*}
\intomo \hthdel^2\eta^2\dx=&-\beta\intomo D^2G_\delta(\nabla\wudel-v)\big(\palpha\nabla\wudel,\nabla\eta^2\otimes (\palpha\wudel-\overline{\palpha\wudel})\big) \dx\\
&+\intomod (\wudel-f)\palpha(\eta^2(\palpha\wudel-\overline{\palpha\wudel}))\dx\\
&+\beta \intomo D^2G(\nabla\wudel-v)(\palpha v,\palpha\nabla\wudel)\eta^2\dx\\
&+\beta \intomo D^2G(\nabla\wudel-v)\big(\palpha v,\nabla\eta^2\otimes (\palpha\wudel-\overline{\palpha\wudel})\big)\dx\\
&=:\widetilde{S}_1+\widetilde{S}_2+\widetilde{S}_3+\widetilde{S}_4.
\end{align*}
We see, that the terms $\widetilde{S}_1$ and $\widetilde{S}_2$ can be treated like the corresponding quantities $S_1$ and $S_2$ above. To $\widetilde{S}_3$ we apply the Cauchy-Schwarz and Young's inequality:
\begin{align*}
|\widetilde{S}_3|&\leq \eps\intomo \wthdel^2\eta^2\dx+c(\eps)\intomo |\nabla v|^2\eta^2\dx\\
&\leq  \eps\intomo \wthdel^2\eta^2\dx+cR^2.
\end{align*}
For $\widetilde{S}_4$ we obtain via similar arguments as for $S_1$ the estimate
\begin{align*}
|\widetilde{S}_4|&\leq \frac{c}{R}\intdb |\nabla v||\nabla \wudel-\overline{\nabla\wudel}|\dx\\
&\leq \frac{c}{R}\left(\intdb |\nabla v|^2\dx\right)^{\frac{1}{2}}\left(\intop_{T}|\nabla\wudel-\overline{\nabla\wudel}|^2\dx\right)^\frac{1}{2}\\
&\leq c\intop_{T}|\nabla^2\wudel|\dx\leq c\intop_{T}\hthdel\widehat{\omega}_\delta\dx.
\end{align*}
Altogether, this suffices to prove the estimate
\begin{align*}
 \intdb \hthdel^2\dx\leq \frac{c}{R}\left[\intop_{T}\thdel^2\dx +R^2\right]^\frac{1}{2}\intop_{T}\hthdel\psidel\dx+cR^\gamma
\end{align*}
for exponents $\gamma\in (0,2)$, which implies
\begin{align*}
\intb \hthdel^2\dx\leq c\frac{1}{\ln\left(\frac{1}{R}\right)^t}\text{ for any $t\geq 1$.}
\end{align*}
From this point on, we can repat the arguments which were used to derive the Hölder continuity of $v$. However, one should note that we have to replace $\sigma_\delta$ with the quantity $\rho_\delta:=DG_\delta(\nabla\wudel-v)$. Then, as in \gr{sig}, we have
\begin{gather*}
\intb|\nabla\rho_\delta|^2\dx=\intb\palpha\rho_\delta\cdot \palpha\rho_\delta\dx\\
=\intb D^2G_\delta(\nabla\wudel-v)\big(\palpha\nabla\wudel-\palpha v,\palpha\rho_\delta\big)\dx\\
\leq c\left[\left(\intb \hthdel ^2 \dx\right)^{\frac{1}{2}}+\left(\intb |\nabla v|^2\dx\right)^\frac{1}{2}\right]\left(\intb |\nabla\rho_\delta| ^2 \dx\right)^{\frac{1}{2}}
\end{gather*}
and therefore
\begin{align*}
\left(\intb|\nabla\rho_\delta|^2\dx\right)^\frac{1}{2}\leq c\left(\intb \hthdel ^2 \dx\right)^{\frac{1}{2}}+cR.
\end{align*}
Since $\lim\limits_{R\rightarrow 0}\ln\left(\frac{1}{R}\right)^\frac{t}{2}R=0$ for any $t\geq 1$, we may neglect the additional term $cR$ and conclude
\[
\left(\intb|\nabla\rho_\delta|^2\dx\right)^\frac{1}{2}\leq c\frac{1}{\ln\left(\frac{1}{R}\right)^\frac{t}{2}}\text{ for any $t\geq 1$.}
\]
 Arguing as for $v$, we thus get an open subset $\widetilde{\Omega_0}\subset\Omega_0$ such that $\Omega_0-\widetilde{\Omega_0}$ has Hausdorff dimension $0$,  and $(u,v)\in C^{1,\alpha}(\widetilde{\Omega_0})\times C^{1,\beta}(\widetilde{\Omega_0},\R^2)$ for every pair $(\alpha,\beta)\in (0,1)\times (0,1)$. The proof of Theorem \ref{Thm1.4} is thereby finished. \qed
\end{section}

\vspace{1cm}

\begin{tabular}{l }
Jan Müller\\
Saarland University\\
Department of Mathematics\\
P.O. Box 15 11 50\\
66041 Saarbrücken\\
Germany\\
jmueller@math.uni-sb.de
\end{tabular}

\end{document}